\documentclass[12pt]{amsart}
\usepackage{amssymb,graphicx,verbatim,epstopdf}

\newtheorem{thm}{Theorem}
\newtheorem{lem}[thm]{Lemma}

\newtheorem{conj}[thm]{Conjecture}
\newtheorem{cor}[thm]{Corollary}
\theoremstyle{definition}

\newtheorem{rem}[thm]{Remark}
\newcommand{\fig}[2]{\includegraphics[scale=#2]{#1.eps}}
\newcommand{\figno}[1]{\centerline{Figure #1}}
\newcommand{\place}[3]
  {\text{\kern#1pt \smash{\raise#2pt\hbox{\rlap{#3}}}\kern-#1pt\kern-1ex}}

\newcommand{\Ker}{\operatorname{Ker}}
\newcommand{\Hom}{\operatorname{Hom}}
\newcommand{\rk}{\operatorname{rank}}

\newcommand{\Z}{{\mathbb{Z}}}
\newcommand{\N}{{\mathbb{N}}}
\newcommand{\Q}{{\mathbb{Q}}}

\newcommand{\R}{{\mathbb{R}}}

\newcommand{\G}{\Gamma}
\newcommand{\La}{\Lambda}
\newcommand{\e}{\varepsilon}
\newcommand{\hra}{\hookrightarrow}

\newcommand{\ra}{\longrightarrow}

\def\smfrac#1#2{\smash{\textstyle\frac{#1}{#2}}}

\def\sm{\setminus}
\def\bn{\begin{enumerate}}

\def\en{\end{enumerate}}
\def\ba{\begin{array}}
\def\ea{\end{array}}
\def\be{\begin{equation}}
\def\ee{\end{equation}}

\numberwithin{thm}{section}

\begin{document}

\title{Non-smoothable four-manifolds with infinite cyclic fundamental group}
\author{Stefan Friedl, Ian Hambleton, Paul Melvin and Peter Teichner}
\address{Universit\'e du Qu\`ebec \`a Montr\'eal, Montr\'eal, Qu\`ebec}
\email{sfriedl@gmail.com}
\address{McMaster University, Hamilton, Ontario}
\email{hambleton@mcmaster.ca}
\address{Bryn Mawr College, Bryn Mawr, Pennsylvania}
\email{pmelvin@brynmawr.edu}
\address{University of California, Berkeley, California}
\email{teichner@math.berkeley.edu}

\thanks{Research partially supported by RTN Network HPRN--
CT--2002--00287: Algebraic K-Theory, Linear Algebraic Groups and Related Structures, NSERC
Discovery Grant A4000, NSF:FRG-0244460, and NSF: DMS-0453957. This joint work started at the Park
City Mathematics Institute 2006, summer program on low-dimensional topology, with support from
IAS}


\begin{abstract}
In \cite{HT}, two of us constructed a closed oriented $4$-dimensional manifold with fundamental
group $\Z$ that does not split off  $S^1\times S^3$. In this note we show that this $4$-manifold,
and various others derived from it, do not admit smooth structures. Moreover, we find an infinite
family of $4$-manifolds with exactly the same properties. As
a corollary, we obtain topologically slice knots that are not smoothly slice in any rational
homology ball.
\end{abstract}

\maketitle

\section{Introduction}

The symmetric matrix
$$
L := \begin{pmatrix}
1+t+t^2 & t+t^2 & 1+t & t \\
t+t^2 & 1+t+t^2 & t & 1+t \\
1+t & t & 2 & 0 \\
t & 1+t & 0 & 2
\end{pmatrix}
$$
has determinant~1 and therefore describes a nonsingular form over $\Z[t]$, see  \cite[p.474]{Q}. Setting $t=x+x^{-1}$ and
$$
\Lambda := \Z[x,x^{-1}] \ \ (=\text{the group ring $\Z[\Z]$}),
$$
$L$ extends to a form on $\Lambda^4$ which is hermitian with respect to the involution $\bar x = x^{-1}$.
In \cite{HT}, the second and fourth authors showed the following:
\begin{thm} \label{thm:old}
$L$ is not extended from the integers.
\end{thm}
Freedman and Quinn \cite{FQ} proved that any nonsingular hermitian form on a finitely generated free $\Lambda$-module can be realized as the intersection form on $\pi_2$ of a closed oriented $4$-manifold with fundamental group $\Z$. Moreover, in the odd case (as for the form $L$ above) there are exactly two such $4$-manifolds realizing a given form, one with nontrivial Kirby-Siebenmann invariant and hence not smoothable.  It follows that there is a unique closed orientable $4$-manifold $M_L$ with $\pi_1M_L=\Z$, intersection form $L$ on $\pi_2M_L$ and trivial Kirby-Siebenmann invariant (so $M_L$ is smoothable after adding copies of $S^2 \times  S^2$). See also Remark~\ref{rem:M_L} for a more concrete construction.

Note that $M_L$ was (before this paper was written) the only known closed orientable topological $4$-manifold with fundamental group $\Z$ and trivial Kirby-Siebenmann invariant  that is not the connected sum of $S^1 \times S^3$ with a simply connected $4$-manifold.  We shall give more examples in Section~\ref{sec:moreforms} below.  Fintushel and Stern constructed a smooth example in \cite{FT}  that is nonsplit in the smooth category but does split off $S^1 \times S^3$ topologically.

It remained an open problem during the last ten years as to whether $M_L$ could be given a smooth structure. In this paper we will prove that $M_L$ is \emph{not} smoothable. Note that Donaldson's Theorem A, \cite{D}, does not apply directly, since the intersection form on $H_2M_L$ is standard.  However, we shall show that it does apply to most cyclic covers, and hence indirectly to $M_L$.

\begin{thm} \label{thm:new}
Let $M_n$ be the n-fold (cyclic) cover of $M_L$. For any $n\geq 1$, the intersection form $L_n$ on $H_2M_n$ is positive definite, odd and of rank~$4n$. If $n=1,2$ this form represents the standard lattice in $\R^{4n}$ but for $n\geq 3$ it is not standard. In particular, none of the $M_n$  admit a smooth structure.
\end{thm}

For $n\geq 3$, the last statement follows directly from \cite{D} since smooth definite $4$-manifolds can only have standard intersection forms. But if $M_1(=M_L)$ or $M_2$ were smooth, then $M_4$ would be smooth as a common covering space.

Theorem~\ref{thm:new} implies Theorem~\ref{thm:old}:  if $L$ was extended from the integers then it would be extended from its augmentation $\e(L)=L_1$, which is the standard form of rank~4. It is easy to see that this implies that all $L_n$ would be standard as well.

Looking closely at the original proof of Theorem~\ref{thm:old}, namely that Lemma~4 in \cite{HT}
holds for forms over the group ring $ \Lambda_n := \Z[\Z/n] $ for all $n>1$, one sees that the
form $L_n$ is not standard as a form with $\Z/n$ action for $n>1$. This means that the
$\Z/n$-equivariant intersection form on $H_2M_n$ is not extended from the integers.  Together with
Theorem~\ref{thm:new}, this gives an interesting algebraic example for $n=2$.

Theorem~\ref{thm:new} suggests the following conjecture, which would reduce the realization
problem for smooth $4$-manifolds with fundamental group $\Z$ to the simply-connected case. As
already conjectured in \cite{HT},  the indefinite case should be purely algebraic and not require
smoothness, whereas the definite case, just as for $M_L$, must use the smooth structure on the
manifold.

\begin{conj} \label{conj:conj}
If $M$ is a closed smooth $4$-manifold with fundamental group $\Z$ then its intersection form on $\pi_2M$ is extended from the integers.
\end{conj}

The proof of Theorem~\ref{thm:new} will be given in Section~\ref{sec:proof}.  In
Section~\ref{sec:forms}, we further analyze the structure of the forms $L_n$, obtaining some
partial results on how they decompose. Section~\ref{sec:knots} contains an application of Theorem~\ref{thm:new} to the study of topologically slice knots that are not smoothly slice:

\begin{cor} \label{cor:knots}
There exists a knot $K$ in $S^3$ such that $0$-surgery on $K$ bounds a smooth $4$-manifold $W$ with fundamental group $\Z$ and such that the intersection form  on $\pi_2W$ is represented by $L$. 
Moreover, any knot with this property has trivial Alexander polynomial
(and is thus topologically slice) but cannot bound a smooth disk in a rational homology
ball.
\end{cor}

 We will prove Corollary~\ref{cor:knots} in Section~\ref{sec:knots}, but here is an explicit example of such a knot $K$.

\bigskip
\centerline{\fig{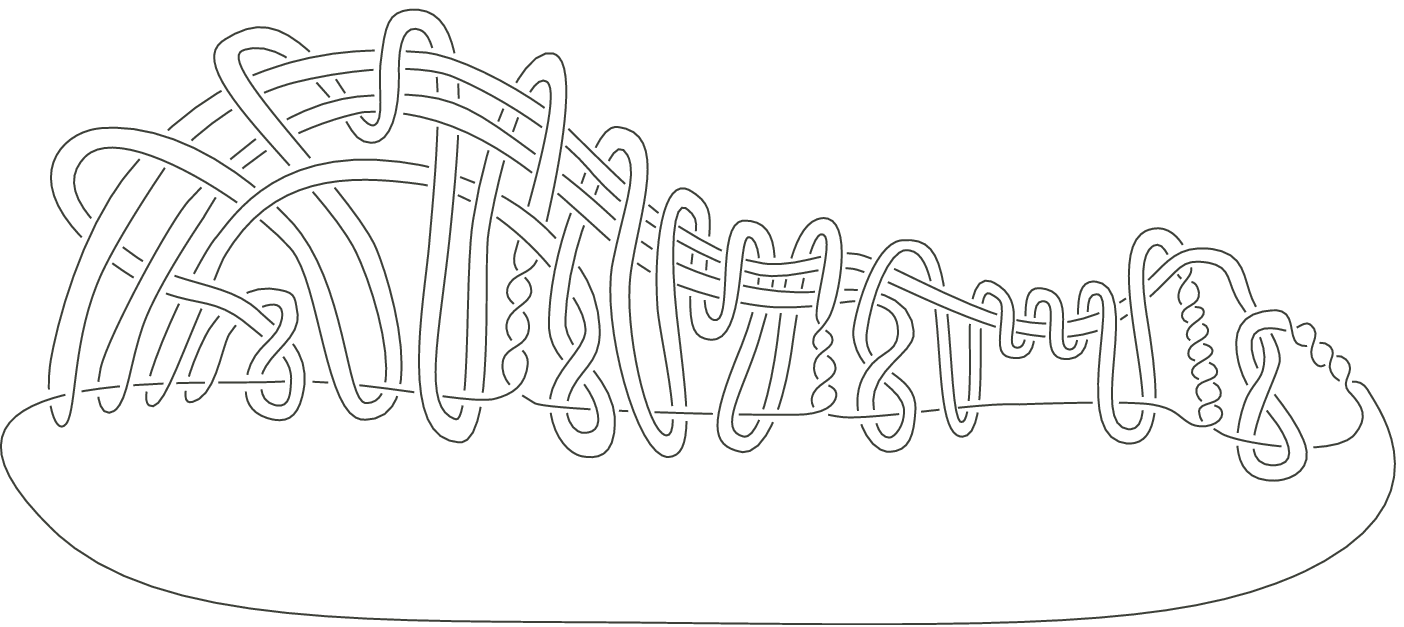}{.8}}
\place{260}{120}{$K$}
\newline
\figno{1: A knot as in Corollary~\ref{cor:knots} }
\medskip

The remaining sections of the paper contain generalizations of our main results. In particular, in
Section~\ref{sec:moreforms} we prove the following theorem. Our Corollary~\ref{cor:knots} continues to hold for all these hermitian forms $L(k)$.

\begin{thm}\label{thm:forms}
There are pairwise non-isomorphic unimodular, odd, hermitian forms $L(k)$, $k\in\N$, on
$\Lambda^4$, with $L(1)=L$ and $\e(L(k))$ standard.  None of these forms is extended from the
integers, and none of the associated closed $4$-manifolds (with infinite cyclic fundamental group
and trivial Kirby-Siebenmann invariant) are smoothable.
\end{thm}

In the final Section~\ref{sec:BS} we show that the phenomenon described in Theorem~\ref{thm:new} can be generalized to other fundamental groups.   In particular, consider a finite aspherical 2-complex with fundamental group $\G$, so
this 2-complex represents a $K(\G,1)$. Given $1\neq g\in \G$, there is an embedding $i_g:\Z\hra\G$ induced by $1\mapsto g$ because $g$ has infinite order (otherwise a nontrivial finite cyclic group would have no higher homology).  We will prove:

\begin{thm}\label{thm:other}
There is a closed oriented $4$-manifold $M$ with $\pi_1(M)=\G$ such that
\begin{enumerate}
\item The intersection form $\lambda$ on $\pi_2M$ (modulo its
radical) is given by extending the form $L$ from $\Z[\Z]$ to
$\Z[\G]$ via $i_g$, but $\lambda$ is not extended from the integers.
\item If there exists an epimorphism $\varphi:\G\to G$ to a finite
group $G$ such that $\varphi(g)$ has order $\geq 3$ and $H_2(\Ker\varphi)=0$, then $M$ does not
admit a smooth structure. In fact the cover $M_\varphi$ of $M$ corresponding to $ \varphi$ has a
non-standard odd, positive definite intersection form.
\end{enumerate}
\end{thm}

Taking $\G=\Z$ we see that Theorem~\ref{thm:new} can be logically seen as a special case of Theorem~\ref{thm:other}. Other classes of groups that satisfy all assumptions of Theorem~\ref{thm:other} include knot groups, fundamental groups of non-orientable surfaces (except for the projective plane but including free groups) as well as almost all Baumslag-Solitar groups, see Lemma~\ref{lem:BS}. We leave it to the reader to formulate (and prove) the appropriate amalgamation of Theorems~\ref{thm:forms} and \ref{thm:other}.

\begin{rem}
The assumption on the vanishing of $H_2(\Ker\varphi)$ in Theorem~\ref{thm:other} is necessary for the last conclusion to hold. In fact, the intersection form on $H_2(M_\varphi)$ contains a metabolic form on $H_2(\Ker\varphi)\oplus H_2(\Ker\varphi)^*$ and hence can only be positive definite under our assumption.
\end{rem}

\textbf{Acknowledgment.} The authors wish to thank the referee for suggesting a better statement and proof of Corollary \ref{cor:knots}.

\section{Proof of Theorem~\ref{thm:new}} \label{sec:proof}

First observe that the signature $\sigma$ and the Euler characteristic $\chi$ are multiplicative in finite covers and that they are both equal to 4 for $M_1=M_L$. But since any finite cover satisfies $\pi_1M_n=\Z$, it actually follows that
\[
\sigma(M_n) = 4n = \chi(M_n) = \rk H_2M_n
\]
and hence the intersection form $L_n$ on $H_2M_n$ is positive definite. The intersection form
$$
L_1= \begin{pmatrix}
7&6&3&2\\
6&7&2&3\\
3&2&2&0\\
2&3&0&2
\end{pmatrix}
$$
on $H_2M_L$ (obtained from $L$ via the augmentation $\e:\Z[x,x^{-1}]\to\Z$, $x\mapsto1$) is odd
and hence $w_2(M_L)\neq 0$ because  $H^2(M_L;\Z/2) \cong \Hom(H_2M_L, \Z/2)$ by the universal
coefficient theorem.  Since $H^2(\Z;\Z/2)=0$ it follows that the second Stiefel-Whitney class is
nontrivial on the universal covering of $M_1$.  Therefore, it must also be nontrivial on $M_n$ and
hence $L_n$ is odd.

As a consequence $L_1$ and $L_2$ are odd definite unimodular forms of ranks~4 and~8, respectively, and hence are standard.  It remains to show that $L_n$ is nonstandard for $n\geq 3$. We shall use the following easy criterion. The converse of this criterion was proven to hold by Elkies \cite{E}. His criterion was used to give a geometrically easier argument for Donaldson's Theorem A via the Seiberg-Witten equations.

\begin{lem}
Let $w\in V$ be a {\em characteristic vector} for a unimodular form on $V$, i.e.
\[
(w,v) \equiv (v,v)\quad \mod 2\quad \forall\ v\in V.
\]
If the form is standard, then $|w|^2 := (w,w) \geq\rk V$.
\end{lem}

\begin{proof}
Let $e_i$ be an orthonormal basis of $V$, i.e.\ $(e_i,e_j)=\delta_{ij}$. Then $w=\sum e_i$ is a
characteristic vector with $|w|^2 = \rk V$ and it suffices to show that any other characteristic
vector has larger norm. But this follows from the inequality
\[
|w + 2v|^2 = |w|^2  + 4(v,w) + 4 |v|^2 \geq |w|^2
\]
which is a consequence of $(v,w) = \sum v_i$ and $|v|^2 = \sum v_i^2$ when $v=\sum v_i\cdot e_i$.
\end{proof}

It therefore suffices to show that there is a characteristic vector $w$ for $L_n$ whose $L_n$-norm satisfies
\[
|w|^2 < 4n = \rk L_n.
\]
We first explain an easy way to think about the forms $L_n$.  As in Section 1, set
$$
\Lambda := \Z[x,x^{-1}] = \Z[\Z] \qquad\text{and}\qquad \Lambda_n := \Lambda/(x^n-1) = \Z[\Z/n].
$$
The natural projection $\Lambda \to \Lambda_n$ converts our matrix $L$ into a matrix $L_{/n}$ over $\Lambda_n$ of rank~4.

\begin{lem}\label{lem:n}
The matrix $L_{/n}$ represents the $\Z/n$-equivariant intersection form on $H_2M_n$ with values in
$\Lambda_n$. The ordinary $\Z$-valued intersection form $L_n$ on $H_2M_n$ is given by the
coefficient of the identity (in the ring $\Lambda_n$) of this form.
\end{lem}

\begin{proof}
Since the higher cohomology groups of $\Z$ vanish, the K\"unneth spectral sequence reduces to the (second) isomorphism
\[
H_2(M_n) \cong H_2(M_1;\Lambda_n) \cong H_2(M_1;\Lambda) \otimes_\Lambda\Lambda_n \cong \pi_2M_1
\otimes_\Lambda \Lambda_n \cong \Lambda_n^{\,4}.
\]
The second statement of the Lemma follows from the well known expression of $\Z[G]$-valued intersections on a $G$-cover $X$
\[
\langle a,b \rangle = \sum_{g\in G} ( \bar g a, b)\cdot   g \quad \forall \ a,b \in H_2X
\]
in terms of ordinary intersection numbers $(\bar g a, b)$  of translates $\bar g a$ with $b$.
\end{proof}

To fix notation, let $V_n$ be the free $\Lambda_n$-module of rank~4 (with underlying free abelian
group of rank~$4n$). Then $L_{/n}$ defines a pairing $ \langle\,,\rangle $ on $V_n$ with values in
$\Lambda_n$ and $L_n$ is the coefficient of the identity element in $\Z/n$ of this pairing. We
write
\[
(v,v') := L_n(v,v') = \langle v, v' \rangle_1 \quad \text{ and hence } \quad |v|^2:= (v,v) =
\langle v, v \rangle_1.
\]
Let $\{e_1, e_2, e_3, e_4\}$ be the $\Lambda_n$-basis for $V_n$ in which we have written the matrix $L$ above. Then one gets a $\Z$-basis $E$ for $V_n$ by multiplying $e_i$ with $x^j$, where $x$ is a generator of $\Z/n$ and $j=1,\dots,n$. One can compute the norms of these basis vectors (noting that $|x^j e_i|^2 = \langle x^j e_i, x^j e_i \rangle_1 = x^j\cdot x^{-j} \langle e_i,e_i \rangle_1 = \langle e_i,e_i \rangle_1 = |e_i|^2$) by looking at the matrix $L$ above.  In particular
$$
|x^j e_i|^2 =|e_i|^2 = \begin{cases}
3,\,5 \text{ or } 7 & \text{for } i=1,2 \\
2 & \text{for } i=3,4.
\end{cases}
$$
The three possibilities in the first case (in which we are computing the coefficient of the identity element in
$\langle e_1, e_1 \rangle = \langle e_2, e_2 \rangle = 1+ t + t^2 = 3 + x + x^{-1} + x^2  + x^{-2} \in \Lambda_n$)
correspond to $n\ge3$, $n=2$ (when $x^2=1$) and  $n=1$ (when $x=1$).   In any case $|e_i|^2$ is always odd for $i=1,2$ and even for $i=3,4$.  This will be used in the lemma below.

Let
$$
N:=1+ x + x^2 +\dots + x^{n-1}
$$
be the norm element in $\Lambda_n$.  It satisfies the identity $N\cdot r$ = $N \cdot \e(r)$ for all $r\in \Lambda_n$, where
$$
\e:\Lambda_n \to  \Z\ , \ \ x\mapsto 1
$$
is the augmentation map. It follows that for any $v,e \in V_n$ we have $\langle N\cdot v, e \rangle = N \langle v,e \rangle = N \e\langle v,e \rangle$ in $\Lambda_n$, and therefore
\[
(N\cdot v, e) =  \langle N\cdot v, e \rangle_1 = ( N \e\langle v,e \rangle)_1 =  \e\langle v,e \rangle \quad\in\Z.
\]

\begin{lem} \label{lem:char}
The vector
\[
w:= N\cdot (e_3 + e_4) \quad \in V_n
\]
is characteristic for $L_n$ and satisfies $|w|^2=4n$.
\end{lem}

\begin{proof}
We need to check that $(w,e) \equiv (e,e) \pmod 2$ for all basis vectors $e =x^j e_i$ in the basis $E$ of $V_n$.   We already calculated the right hand side, and for the left hand side we have
\[
(w,x^j e_i) = \e\langle e_3 + e_4, e_i \rangle =
\begin{cases}
5 & \text{for } i=1,2 \\
2 & \text{for } i=3,4
\end{cases}
\]
by referring to the intersection matrix $L_1=\e(L)$ displayed above.  These numbers are indeed odd for $i=1,2$ and even for $i=3,4$, as required.  Finally, we calculate the norm
\[
|w|^2 =(w,w) = \e(N \langle e_3+e_4,e_3+e_4 \rangle)= \e(N) \e(2+2) = 4n.
\]
\vskip -.25in
\end{proof}

\smallskip
Now consider the characteristic vector $w_1:= w-2e_1$ and calculate
\[
|w_1|^2 = |w|^2  - 4\left((w,e_1) - |e_1|^2\right)= 4n - 4\left(5- |e_1|^2\right).
\]
For $n=2$ we have $|e_1|^2 = 5$, and so $|w_1|^2=4n=8$, consistent with $L_2$ being standard.  But for $n\geq 3$ we have $|e_1|^2 = 3$, and hence $|w_1|^2=4n-8$, which means that $L_n$ cannot be standard.

\section{Intersection forms} \label{sec:forms}

In the previous section we showed that the intersection forms $L_n$ on
$V_n= H_2(M_n;\Z)$ ($\cong \Z^{4n}$) are not standard for any $n\ge3$,
i.e.\ not equivalent to the diagonal form $I_{4n}$.   This section
initiates a further investigation of the structure of these forms, and
thus implicitly the topological structure of the $4$-manifolds $M_n$.

In what follows, a free abelian group of finite rank with a positive
definite unimodular form will be referred to simply as a {\it
lattice}.   By a classical theorem of Eichler (see \cite[6.4]{MH}), any
lattice splits uniquely as an orthogonal sum of indecomposable
lattices.   There is a unique indecomposable lattice in each rank $1$,
$8$, $12$, $14$ and $15$, but after that the number grows dramatically
with the rank.

Of special interest here are the lattices $\G_{4m}$ (of rank $4m$) which
are indecomposable for $m\ge3$.  In particular $\G_8$ (also known as
$E_8$) and $\G_{12}$ (also known as $D_{12}^+$ as in Conway and Sloan's
treatise \cite{CS}) are the smallest non-trivial indecomposable
lattices.  Explicitly, $\G_{4m}$ consists of all vectors in $\R^{4m}$
whose coordinates are all integers or all half-integers with an even
integer sum.  It is spanned by the vectors $v_i+v_j$ and
$\frac12(v_1+\cdots+v_{4m})$ (for an orthonormal basis $v_1,\dots,v_{4m}$).

Before stating any results, we recall some standard notions (and introduce some new ones)
concerning a lattice $V$ with inner product $(\,, )$.  The {\it norm} of $v\in V$ is the inner
product $(v,v)$, also denoted $|v|^2$.  A vector $w\in V$ is {\it characteristic} if
$(w,v)\equiv|v|^2 \mod 2$ for all $v\in V$.  Among the characteristic vectors, those with the
smallest norm will be called {\it minimal}.  The total number $\mu(V)$ of minimal characteristic
vectors in $V$ is a useful invariant of $V$.  It clearly multiplies under orthogonal sums, since a
vector $(u,v)\in U\oplus V$ is (minimal) characteristic if and only if $u\in U$ and $v\in V$ are
(minimal) characteristic.  For example  we have $\mu(I_k)=2^k$ (where $I_k$ denotes the standard
lattice) and $\mu(\G_8)=1$ (since $\G_8$ is even).

Define the {\it defect} of the lattice to be
$$
d(V) = \textstyle\frac18(\rk(V)-|w|^2)
$$
for any minimal characteristic vector $w$ in $V$.   It is an integer, by
a classical lemma of van der Blij (see \cite[5.2]{MH}), and clearly adds
under orthogonal sums.  For example $d(\G_{4m}) = \lfloor m/2 \rfloor$.
Indeed for odd $m$ the minimal characteristic vectors are exactly the
vectors $\pm2v_i$ (and so $\mu(\G_{4m}) = 8m$) while $0$ is the unique
minimal characteristic vector when $m$ is even.

Elkies shows in \cite{E}  that the defect $d(V)$ is always nonnegative,
and equal to $0$ if and only if $V$ is standard.  It is therefore a
measure of non-diagonalizability.  In \cite{E'}, he goes on to prove
that there are only finitely many indecomposable lattices of defect one
(including $\G_8$ and $\G_{12}$) and raises the question of whether this
is also the case for larger defects.  (This has been confirmed for
defects $2$ and $3$ by Gaulter \cite{G}.)

\begin{thm} \label{thm:defect}
The defect of the lattice $V_n \ (= H_2(M_n;\Z)$ with the form $L_n)$
satisfies the inequality $\lfloor n/3 \rfloor \le d(V_n) < n/2.$
\end{thm}

\begin{proof}
The result is obvious for $n=1$ or $2$, so assume $n\ge3$.  We use the
notation from the previous section, in particular, $e_1,\dots,e_4$ is a
basis for which the intersection form $L$ is given by the matrix on page
1 and $w = N\cdot (e_3+ e_4)$. We consider the characteristic vector
$$
w_0 := w-2(1+x^3+\cdots+x^{3(\lfloor n/3 \rfloor - 1)})e_1
$$
in $V_n$.  This vector is of the general form
$$
w_a := w-2a(x)e_1
$$
where $a(x) = a_0 + a_1x + \cdots + a_k x^k$ is an integer polynomial.
The norm of $w_a$ can be written in terms of the quadratic expressions
$a^i := a_0a_i+a_1a_{i+1}+\cdots+a_ka_{i-1}$ (i.e.\ the dot product of
the coefficient vector $(a_0,a_1,\dots,a_k)$ with its $i$th shift
$(a_i,a_{i+1},\dots,a_{i-1}$)) as follows:
$$
\begin{aligned}
|w_a|^2 &= 4n - 4\left(\e\langle e_3+e_4,a(x)e_1\rangle -
\left(a(x)a(x^{-1})\langle e_1,e_1\rangle\right)_1\right) \\
&= 4n-4\left(5a(1) - 3a^0 - 2(a^1+a^2)\right).
\end{aligned}
$$
For $w_0$, the coefficient vector consists of a sequence of $\lfloor n/3
\rfloor$ ones separated by and terminating with at least two zeros, and
so $a(1)=a^0=\lfloor n/3 \rfloor$ and $a^1=a^2=0$.  Thus $|w_0|^2 =
4n-8\lfloor n/3\rfloor$, and so by definition $d(V_n)\ge\lfloor n/3\rfloor$.

To obtain the upper bound, note that $d(V_n) \le \frac18\rk(V) = n/2$,
and the inequality is strict because the form on $V_n$ is odd.
\end{proof}

\begin{cor} \label{cor:largen}
If $V_n$ splits off an $I_k$ summand for some $k$, then $k\le
4n-8\lfloor n/3\rfloor$.
\end{cor}

This is immediate from the theorem and the additivity of the defect.
Observe that the bound $4n-8\lfloor n/3\rfloor$ is of the order of
$4n/3=\rk(V_n)/3$ for large $n$, and so this result roughly states that
at least two-thirds of the form $V_n$ cannot be diagonalized.

Unfortunately this gives no new information for $n<6$ beyond that
contained in Theorem~\ref{thm:new}.  Special arguments can be used,
however, to further restrict the forms for small $n$.  Indeed for $n=3$
or $4$ they can be identified precisely, using Conway and Sloan's root-system
labeling $($also see Conway's lovely little book
\cite[pp.\,53--58]{C}$)$.

\begin{thm} \label{thm:smalln}
$V_3$ and $V_4$ are the unique odd indecomposable lattices in dimensions
$12$ and $16$.  In particular,  we have $V_3\cong\G_{12}$ and $V_4\cong
D_8^2[(12)]$.\end{thm}

\begin{proof}
The lattices of rank $12$ are $\G_{12}$, $\G_8\oplus I_4$ and $I_{12}$,
and so by Theorem~\ref{thm:new} it suffices to eliminate $\G_8\oplus
I_4$.  But
$$
\mu(\G_8\oplus I_4) = 1\cdot2^4 = 16 \quad\text{while}\quad \mu(\G_{12})=24
$$
and so one need only produce $17$ or more minimal characteristic vectors
in $V_3$.  In fact the complete list is $\pm(w - 2x^i e)$ for $i=0,1,2$
and $e = e_1$, $e_2$, $e_1+e_4$ or $e_2+e_3$, as is readily verified.
(An explicit isomorphism $V_3\cong \G_{12}$ is easily deduced from this.)

For $V_4$ we take a different tact, focusing on vectors of norm $2$.  A
theorem of Witt asserts that the norm 2 vectors in any lattice span a
sub-lattice isomorphic to a direct sum of root lattices $A_n$, $D_n$ and
$E_n$.  From the classification in \cite[\S16.4]{CS} of lattices of rank
16,  it suffices to find a copy of $D_8\oplus D_8$ in $V_4$, where $D_8$
is the lattice whose intersection matrix with respect to a suitable
basis $v_1\dots,v_8$ of norm 2 vectors is given by the corresponding
Dynkin diagram:
\newsavebox{\diagram}
\savebox{\diagram}{\begin{picture}(0,0)
\put(-72,14){\circle*{4}}
\put(-48,14){\circle*{4}}
\put(-24,14){\circle*{4}}
\put(0,14){\circle*{4}}
\put(24,14){\circle*{4}}
\put(48,14){\circle*{4}}
\multiput(63,-2)(0,32){2}{\circle*{4}}
\put(-72,14){\line(1,0){120}}
\put(49,16){\line(1,1){14}}
\put(49,12){\line(1,-1){14}}
\put(-75,2){$v_1$}
\put(-51,2){$v_2$}
\put(-27,2){$v_3$}
\put(-3,2){$v_4$}
\put(21,2){$v_5$}
\put(41,2){$v_6$}
\put(68,-4){$v_7$}
\put(68,28){$v_8$}
\end{picture}}
\bigskip
\begin{center}
\begin{picture}(20,30)
 \setlength{\unitlength}{1.5pt}
 \put(0,0){\usebox{\diagram}}
\end{picture}
\end{center}
\bigskip
Here two nodes are joined by an edge or not according to whether their
associated vectors have inner product $1$ or $0$.  One easily checks
that in $V_4$ the vectors
$$
v_1= x^2e_4, \quad
v_2 = x^2(-e_1+e_2+e_3), \quad
v_3 = x^2e_3, \quad
v_4 = x^2e_1-e_2, \quad
v_5 = e_1-e_2
$$
\vskip-.3in
$$
v_6 = e_3, \quad
v_7 = -e_1+e_2+e_3-e_4, \quad
v_8 = w-(e_1+e_2+x^2(e_3+e_4))
$$
span a copy of $D_8$, and the vectors $w_i=x v_i$ (for $i=1,\dots,8$)
span an orthogonal copy of $D_8$.
\end{proof}

\section{Alexander polynomial one knots that are not slice} \label{sec:knots}

In this section we shall prove Corollary~\ref{cor:knots}, starting with the second part. Let $K$ be a knot  in $S^3$ such that the $3$-manifold $N$ obtained by $0$-surgery on $K$ bounds a smooth $4$-manifold $W$ with fundamental group $\Z$ and such that the intersection form  on $\pi_2W$ is represented by $L$. 

Since the intersection form on $H_2(W;\Lambda) \cong\pi_2W$ is non--singular it follows immediately that $\Delta_K(t)=1$, hence $K$ is topologically slice. Moreover, $H_2(W) \cong \pi_2W \otimes_\Lambda \Z$ is torsionfree and therefore the same is true for $H_1(W,N) \cong H^3(W)$. This, together with the non--singularity of $L\otimes_\Lambda \Z$  implies that the induced map $H_1(N)\to H_1(W)$ is an isomorphism.

Now assume that  $K$ bounds a smooth disk $D$ in a rational homology ball $B$. Then let $C=B\sm
\nu D$ and $X=S^3\sm \nu K$. Note that $X$ and $C$ have the rational homology of a circle.
Consider the \emph{smooth} 4--manifold
\[
M:= W\cup_\partial C.
\]
We have inclusion induced isomorphisms \[ \Z\cong H_1(X)\xrightarrow{\cong}
H_1(N)\xrightarrow{\cong} H_1(W)\xrightarrow{\cong} H_1(M)/\{\mbox{torsion}\}\xleftarrow{\cong}
H_1(C)/\{\mbox{torsion}\}.\] Given $n\in \N$ we denote the covers corresponding to the
homomorphism $\Z\to \Z/n$ by the subscript $n$. Now let $n$ be a prime power larger than two such
that $H_*(C)$ has no $n$--torsion. It follows immediately from \cite[Proof~of~Lemma~2.3]{F} (cf.
also \cite[p.~184]{CG}) that $X_n$ and $C_n$ are still rational homology circles and that the
projection maps $X_n\to X$ and $C_n\to C$ give isomorphisms of rational homology.

It is well--known that $H_1(N_n)\cong  H_1(X_n)$. By Poincar\'e duality it now follows that $N_n$
is a rational homology $S^1\times S^2$. Using $H_1(N)\xrightarrow{\cong} H_1(W)$ it follows from
the above discussion and from the Meyer--Vietoris sequence corresponding to $M_n=W_n\cup_{N_n}C_n$
that the inclusion maps induce an isomorphism
\[ H_2(W_n)/\{\mbox{torsion}\}\xrightarrow{\cong} H_2(M_n)/\{\mbox{torsion}\}.\]
Since $\pi_1(W)=\Z$ we can apply Lemma \ref{lem:n} to conclude that the intersection form on $W_n$
(and hence on the smooth 4--manifold $M_n$) is given by $L_n$. But
since $L_n$ is positive and non--standard by Theorem \ref{thm:new} we can again use Donaldson's
theorem A to see that the $M_n$ cannot exist smoothly and hence $K$ cannot be smoothly slice to
begin with. This concludes the proof of the second part of Corollary~\ref{cor:knots}.

To prove the first part of Corollary~\ref{cor:knots} it suffices to construct a smooth
$4$-manifold $W$ with the required properties and check that $\partial W$ is 0-surgery on a knot.
Start with one 0- and one 1-handle (which we draw as a dotted unknot $U$ in $S^3$). Then $W_J$ is
obtained by attaching four 2-handles along a framed link $J$ in $S^1 \times D^2$, the complement
of $U$. The link $J$ is chosen so that
\begin{enumerate}
\item $J$ is the unlink if one ignores $U$,
\item $J$ has trivial linking numbers with $U$ and
\item $J$ represents the matrix $L'$.
\end{enumerate}
Here $L'$ is obtained from $L$ by row and column operations that have the effect of changing
$\e(L)$ to the identity matrix $\e(L')$. Note that property (2) implies that $W_J$ is homotopy
equivalent to $S^1 \vee 4 \  S^2$ and a $\Z[\Z]$--basis of $\pi_2W_J$ is given by the cores of the
four 2-handles together with null-homotopies of $J$ in the complement of $U$. These four
null-homotopies read off an intersection matrix depending on the crossing changes necessary.
Property (3) requires that this matrix is given by $L'$, in particular $J$ is 1-framed (when
ignoring $U$).

It is easy to see that any matrix can be realized by a link $J$ with properties (1) and (2). The
last step is to prove that $\partial W_J$ is 0-surgery on a knot: This uses property (1) and the
fact that $J$ is 1-framed. As a consequence one can blow down $J$ completely without changing the
boundary! If one draws a picture with four disjointly embedded disks bounding $J$ and punctured by
$U$, then the blow down procedure puts a full twist into all strands that go through each of these
four disks. Therefore, $U$ turns into a knot $K$ (and the framing has turned from a dot to a zero
when studying the boundary alone). This is our non-smoothly slice knot with  trivial Alexander
polynomial, in fact, it's a whole family of such knots.

One example of such a knot $K$ is shown in Figure 1. It is constructed from the unknot by four
finger moves.  The self linking and twisting of the fingers produce the diagonal entries of $L$,
while the linking between fingers produce the off diagonal entries. To see this, blow up four $+1$
curves (the link $J$ in the discussion above) to unhook each of the finger tips.  This transforms
$K$ into an unknot $U$.  Now pull $J$ back along the fingers so that $U$ appears as the round
unknot while the components of $J$ follow the original fingers, clasping at both ends.

To calculate the associated intersection form, first introduce the notation
$$
q_n = 2 x^n Ð (x^{n+1} + x^{n-1}) \quad\text{and}\quad c_n = n q_0 Ð (q_1 + q_{-1})
$$
so the matrix $L$ can be written as the sum of the identity matrix with a matrix in $2\times2$ block form in which each block has all entries equal, namely to $c_1$ and $c_3$ for the two diagonal blocks, and $- c_2$ for both off diagonal blocks.

Now orient $U$ and all the components of $J$ counterclockwise and label each arc of $J$ with a
power of $x$ in a natural way:  Start with $1$'s right before the clasps, and then proceed along
$J$ (following the orientation) multiplying by $x$ (or $x^{-1}$) each time $J$ links $U$
positively (or negatively).  Locally, each finger then has two oppositely oriented strands that
are labeled by adjacent powers of $x$, say $x^m$ and $x^{m-1}$, which we abbreviate by simply
labeling the finger at that point with the higher exponent $m$.  Also any crossing between fingers
can be given a sign, namely the sign of the crossing between the strands with the higher exponents
of $x$. With these conventions, a positive finger crossing in which the $i$--th finger with label
$m$ passes over the $j$th with labels $n$ will contribute $q_{m-n}$ to the $ij$th entry of $L$,
and an analogous negative crossing contributes $-q_{m-n}$. Self crossings and twisting of the
$i$th finger contribute in a similar way to the diagonal entries. The calculation is then
straightforward.

\begin{rem}\label{rem:M_L}
Note that any knot $K$ whose 0-surgery $N$ bounds a $4$-manifold $W$ as in
Corollary~\ref{cor:knots} has trivial Alexander polynomial. Therefore, $K$ is topologically
$\Z$-slice and hence there is a topological $4$-manifold $C$ that is a homotopy circle and has
boundary $N$. To construct our non-smoothable 4-manifold in the title slightly more concretely, we
only have to prove that
\[
M_L:= W_J\cup_\partial C
\]
has trivial Kirby-Siebenmann. This follows from the additivity of the (relative) Kirby-Siebenmann
and the fact that $K$ has trivial Arf invariant: $\partial C$ also bounds a spin manifold $W'$
with signature zero and
\[
KS(W\cup_\partial C) = KS(C,\partial C) = KS(W'\cup_\partial C) = 0
\]
because on spin manifolds the Kirby-Siebenmann invariant equals the signature divided by~8.
\end{rem}

\section{More forms} \label{sec:moreforms}

Recall that our basic form $L:\Lambda^4\times\Lambda^4\to\Lambda$ (where $\Lambda=\Z[x,x^{-1}]$) was obtained by substituting $x+x^{-1}$ for $t$ in the matrix displayed at the beginning of this paper.  In this section we investigate the unimodular forms
$$
L(a) := \begin{pmatrix}
1+a+a^2 & a+a^2 & 1+a & a \\
a+a^2 & 1+a+a^2 & a & 1+a \\
1+a & a & 2 & 0 \\
a & 1+a & 0 & 2
\end{pmatrix}
$$
obtained by substituting other elements $a\in \Lambda$ for $t$.  (We shall always use the same notation for a matrix and its associated form.)   Note that $L(a)$ is hermitian with respect to the involution $x\mapsto x^{-1}$ if and only if $a=\bar a$, where $\bar a$ denotes the image of $a$ under this involution (called the ``conjugate" of $a$).   So we assume that $a=\bar a$ and can formulate the precise version of Theorem~\ref{thm:forms} in the introduction.

\begin{thm}\label{thm:old'}
Define inductively $b_1=1$ and $b_{k+1}=4b_k+1$, and set
\[
L(k) := L(x^{b_k}+x^{-b_k}).
\]
Then the forms $L(k)$ for $k=1,2,\dots$ are pairwise non-isomorphic.  None of these forms is extended from the integers, and none of their associated $4$-manifolds are smoothable.
\end{thm}

As with $L=L(1)$ the proof requires a study of some related forms.  First some notation.  For any
$a=\sum a_ix^i\in\Lambda$, recall that $\e(a) = \sum a_i \in \Z$ and $(a)_1=a_0$.  Thus
$\e:\Lambda\to\Z$ is the usual augmentation map, and $(a)\mapsto (a)_1$ defines a linear
(non-multiplicative) projection $\pi:\Lambda\to\Z$. We consider also the corresponding maps
$\Lambda_n\to \Z$, denoted by the same names for any $n\in \N$, where $\Lambda_n =
\Lambda/(x^n-1)$.

Now for any (self-conjugate) $a\in\Lambda$ and any $n$, the matrix $L(a)$ can be viewed as a matrix over $\Lambda_n$ in which case it  is denoted by $L(a)_{/n}$, as is the associated form $\Lambda_n^4\times\Lambda_n^4\to\Lambda_n$.  These forms in turn induce integral forms $L_n(a):\Lambda_n^4\times\Lambda_n^4\to\Z$ of rank $4n$ (since $\Lambda_n^4\cong\Z^{4n}$) by composition with $\pi$, that is $L_n(a)(v,w) = (L(a)_{/n}(v,w))_1$.

Note that $L(a)_{/1}=L_1(a) = \e(L(a)) = \e(L(a)_{/n})$ for all $n$.

\begin{lem}\label{lem:general}
For any $a\in\Lambda$ with $\bar a = a$, the forms $L_n(a)$ are positive definite, odd and of rank $4n$.  In particular the forms $L_1(a)$ and $L_2(a)$ are standard.  Furthermore, if $L(a)_{/n}$ is extended from the integers, then $L_n(a)$ is standard.
\end{lem}

\begin{proof}
We first show that $L_1(a)$ is standard. The $(1,1)$-entry of its defining matrix is
$1+\e(a)+\e(a)^2$ which is clearly odd for any $a$.  By the classification of odd forms of rank 4
it follows that $L_1(a)$ is standard if and only if it is positive definite.  But this can be
checked over $\Q$.  First observe (following Quebbemann's definition \cite[p.~474]{Q}) that
$$
L(a) = \begin{pmatrix} \frac12(I+B(a)^2) & B(a) \\ B(a) & 2I \end{pmatrix}
\qquad
\text{where}
\qquad
B(a) := \begin{pmatrix} 1+a & a \\ a & 1+a \end{pmatrix}.
$$
Now a quick calculation gives
$$
PL(a)\bar P^t = D \qquad\text{where}\qquad P = \begin{pmatrix}I & - \frac12 B(a) \\ 0 & I \end{pmatrix} \quad\text{and}\quad D = \begin{pmatrix} \frac12 I & 0 \\ 0 & 2I \end{pmatrix}
$$
and so applying the augmentation map $\e:\Q[x,x^{-1}]\to\Q$ to this equation we see that $L_1(a)$
is positive definite over $\Q$.  (This also shows that $L(a)$ considered as a form on
$(\Lambda\otimes_{\Z}\Z[\frac12])^4$ is extended from a form defined on $\Z[\frac12]^4$.)  The
argument given at the beginning of Section~\ref{sec:proof} now yields all but the last statement
of the lemma.

Now assume that $L(a)_{/n}$ is extended from the integers.  Then it must be extended from its augmentation $\e(L(a)_{/n}$),
which we have just seen is standard.  Hence $L_n(a)$ is standard.
\end{proof}

\begin{lem}\label{lem:specific}
Let
$$
a=a_0+\sum_{\ell=1}^m a_\ell(x^\ell+x^{-\ell}).
$$
If $a_0^2+2\sum a_\ell^2 < a_0 + 4\sum a_\ell$ (where all sums are from $1$ to $m$) then $L_n(a)$ is not standard for any $n>4m$.
\end{lem}

Note that this lemma applies to all elements $a$ for which all $a_i\in\{0,1\}$ and $a_\ell\ne0$
for at least one $\ell>0$, and in particular to $a = x^k+x^{-k}$ for $k>0$.

\begin{proof}
First observe that all exponents in powers of $x$ appearing in the matrix $L(a)$ are between $-2m$
and $2m$, and so the condition $n>4m$ precludes any cancelation when passing to the quotient
$L(a)_{/n}$.  For notational convenience write $(u,v) = L_n(a)(u,v)$ for the inner product and
$|v|^2=(v,v)$ for the associated norm.

We now proceed as in the proof of Theorem~\ref{thm:new}.  Let $N=1+x+\cdots + x^{n-1}$ and $w=N(e_3+e_4)$.   We will show that $w_1=w-2e_1$ is characteristic of norm $<4n$.

The argument that $w$ is characteristic, which of course implies that $w_1$ is as well, is exactly the same as the argument in Section~\ref{sec:proof}.  We need only observe that for any $j=0,\dots,n-1$, the inner products $(x^je_i,x^je_i)$ are odd for $i=1,2$ and even for $i=3,4$ \ (in particular equal to
$$
\textstyle (1+a+a^2)_1 = 1+a_0+a_0^2+2\sum a_\ell^2
$$
in the first case, using the condition $n>4m$, and $(2)_1=2$ in the second) \ and that the same is true of the inner products $(w,x^je^i) = \e(L(a)_{/n}(e_3+e_4,e_i))$ (which equal $$
\textstyle \e(1+2a)= 1+2\e(a) = 1+2(a_0+2\sum a_\ell)
$$
in the first case and $\e(2)=2$ in the second).

The norm $|w_1|^2 = |w|^2-4((w,e_1)-|e_1|^2)$ which equals
$$
\textstyle 4n-4\left((1+a_0+a_0+4\sum a_\ell)-(1+a_0+a_0^2+2\sum a_\ell^2)  \right)
$$
by the calculations above.  The condition on $a$ is exactly what is needed to show that this is less than $4n$.
\end{proof}

To prove Theorem~\ref{thm:old'}, we first show that the form $L(k) = L(x^{b_k}+x^{-b_k})$ is not isomorphic to any $L(j)$ for $j<k$.  Indeed, Lemma~\ref{lem:specific} shows on the one hand that the form $L_n(j)$ is non-standard for $n>4b_j$, which implies by Lemma~\ref{lem:general} that $L(j)_{/b_k}$ is not extended from the integers since $b_k >4b_j$.  On the other hand the matrix $L(k)$ has entries in $\Z[x^{b_k},x^{-b_k}]$, whence $L(k)_{/b_k}$ has integral entries and so is extended from the integers.

Next we observe that $L(k)$ is not extended from the integers.  For if it were, then each $L(k)_{/n}$ would be as well, which would imply by Lemma~\ref{lem:general} that all the forms $L_n(k)$ would be standard, contradicting Lemma~\ref{lem:specific}.

The final statement in the theorem follows from Lemma~\ref{lem:specific}, exactly as in the proof of Theorem~\ref{thm:new}. \qed

\begin{rem}
The argument in the previous paragraph shows more generally that if $L:\Lambda^n\times\Lambda^n\to\Lambda$ is any form extended from the integers with $L_1$ standard, then all the integral forms $L_n$ are standard.  It is conceivable that the converse holds as well.  This would imply Conjecture~\ref{conj:conj} in the definite case.
\end{rem}

\section{More general fundamental groups} \label{sec:BS}

In this section we want to prove Theorem~\ref{thm:other} from the introduction. But first we recall a class of groups to which the theorem actually applies.
Let $k\in \Z$ and consider the Baumslag--Solitar groups
\[
\G_k:=\langle a,b | aba^{-1}=b^k\rangle.
\]
Note that $\G_0=\Z, \G_1=\Z^2$ and that for $k\neq 0$ we have a semi-direct product decomposition (where $a$ generates the quotient $\Z$ and $b$ corresponds to $\frac{1}{k}$):
\[
\G_k\cong \Z\left[\smfrac{1}{k}\right] \rtimes \Z.
\]
Here $n\in \Z$ acts on $\Z[\frac{1}{k}]$ by multiplication by
$\left(\frac{1}{k}\right)^n$.

\begin{lem}\label{lem:BS}
The 2-complex corresponding to the above presentation of $\G_k$ is aspherical. Let $N_n$ be the
index $n$ normal subgroup of $\G_k$ generated by $a^n$ and $b$. If $k\neq 1$, and $n$ is odd in
the case $k=-1$, then for all other $k,n$ we have
\[
H_i(N_n) = 0 \quad \forall \ i>1.
\]
\end{lem}

\begin{proof}
The presentation of $\G$ has a unique relation which is not a proper power. Then the corresponding
2-complex is aspherical \cite{Ly}, just like for surface groups. It is easy to see that $N_n \cong
\Z[\frac{1}{k}] \rtimes n\cdot\Z$. Then the Wang sequence shows the homology result.
\end{proof}

Before we prove Theorem~\ref{thm:other}, we collect some useful information that applies to closed
oriented $4$--manifolds $M$ with arbitrary fundamental group $\G$. Denote by $\La$ the group ring
$\Z[\G]$ and let $A$ be a $\La$-algebra.

\begin{lem}\label{lem:rad}
Consider the equivariant intersection form with coefficients in $A$
\[
\lambda_A : H_2(M;A) \ra \Hom_A(H_2(M;A), A).
\]
given by Poincar\'e duality composed with the Kronecker evaluation.
If $A=\La$ then the radical $\Ker(\lambda_{\La})$ is isomorphic to $H^2(\G;\La)$. If $A=\Z[G]$ where $G=\G/N$ is a quotient group, then the image of $H^2(\G;A)$ in $H^2(M;A)\cong H_2(M;A)$ is contained in the radical $\Ker(\lambda_A)$ if $H_2(N)$ is finite.
\end{lem}

\begin{proof}
The universal coefficient spectral sequence gives an exact sequence
\[
0 \ra H^2(\G, \La) \ra H^2(M;  \La) \ra \Hom_{ \La}(H_2(M; \La),  \La)
\]
Applying Poincar\'e duality $H^2(M;\La)\cong H_2(M;\La)$ and identifying the composition with
$\lambda$ (respectively $\lambda_{\G/N}$) yields the first result. For the second statement we
note that by naturality of the evaluation map, it suffices to prove the vanishing of
\[
\Hom_A(H_2(\G;A), A) = \Hom_A(H_2(N;\Z), A).
\]
This follows from our assumption that $H_2(N;\Z)$ is finite.
\end{proof}

\noindent {\em Proof of Theorem~\ref{thm:other}.} The construction of the $4$-manifold $M$ is very
similar to that of Section~\ref{sec:knots}. Start with a $4$-dimensional thickening $K$ of the
finite aspherical 2-complex, i.e. $K$ is a handlebody with handles of index~$\leq 2$ corresponding
to the cells of the 2-complex. We may assume that there is a 1-handle $h_g$ that corresponds to
the given element $1\neq g\in\G$.

Next we attach four more 2-handles in a neighborhood of the meridian to the 1-handle $h_g$. This is done just like in Section~\ref{sec:knots} so that the intersections between the null homotopies for the handles represent the matrix $L$. In particular, these 4-handles are attached homotopically trivially, so that the resulting $4$-manifold $X$ has the homotopy type
\[
X\simeq K \vee 4\cdot S^2.
\]
If we had attached these four 2-handles to $S^1 \times D^3$, we showed in Section~\ref{sec:knots} that the resulting $4$-manifold has boundary equal to 0-surgery on a knot with trivial Alexander polynomial. By Freedman's theorem, this knot is $\Z$-slice, i.e. the 0-surgery also bounds a topological $4$-manifold $C$ which is a homotopy circle. Since we have attached the four 2-handles in $X$ only to a meridian of $h_g$, it follows that the boundary $\partial X$ also bounds a $4$-manifold $Y$ that is obtained from $C$ by adding 1- and 2-handles (corresponding to the handles in $K$). Since $C$ is a homotopy circle, it follows that $Y$ is another $K(\G,1)$. We define
\[
M := X\cup_{\partial} Y
\]
and we claim that it has the two properties stated in
Theorem~\ref{thm:other}. The inclusion $Y\hra M$ induces an exact
sequence (with coefficients in any $\La$-module $A$)
\begin{equation}\tag{$*$} \label{equ1}
H_2(Y;A) \ra H_2(M;A) \overset{j}{\ra} H_2(M,Y;A) \overset{0}{\ra} H_1(Y;A) \overset{\cong}{\ra}
H_1(M;A).
\end{equation}
The rightmost map is an isomorphism since the inclusion $Y\hra M$
induces an isomorphism on fundamental groups.

We now turn to the proof of part (1) of Theorem~\ref{thm:other}. Consider the exact sequence~(\ref{equ1}) with  $A=\La=\Z[\G]$. Note that $H_2(Y;\La)=0$ since $Y$ is a $K(\G,1)$. Therefore, $j$ is an isomorphism and we get by excision and
Poincar\'e duality
\[
H_2(M;A) \cong H_2(M,Y;A) \cong H_2(X,\partial X;A) \cong H^2(X;A)
\cong H^2(K;A) \oplus H^2(S^2;A)^4
\]
where $ H^2(K;\La)\cong H^2(\G;\La) $ is the radical of the $\La$-valued intersection form on
$\pi_2M$ by Lemma~\ref{lem:rad}. The intersection form on $H^2(S^2;\La)^4$ (and hence the
intersection form $\lambda_\La$ on $\pi_2M$ modulo its radical) is by construction given by
extending the form $L$ from $\Z[\Z]$ to $\Z[\G]$ via $i_g$. Since $g$ has infinite order, the
proof of \cite[Lemma~4]{HT} can easily be seen to carry through to show that $\lambda_\La$ is not
extended from the integers. Alternatively the fact, proven in (2), that $\lambda_\La$ gives rise
to a non--standard form over $\Z$ also shows that $\lambda_\La$ is not extended.

For part (2) of Theorem~\ref{thm:other}, let  $A =\Z[G]=\Z[\G/N]$ and recall that by assumption
\[
H_2(Y;A) \cong H_2(\G;A) \cong H_2(N;\Z)=0.
\]
By Lemma~\ref{lem:rad}, the group $H^2(K;A) \cong H^2(\G;A)$ lies in the radical of $\lambda_A$
and so the above exact sequence~(\ref{equ1}) and the exact sequence slightly below it show that
the $A$--valued intersection pairing $\lambda_A$ (modulo its radical) is given by extending the
form $L$ from $\Z[\Z]$ to $A=\Z[G]$ via the map $\Z\to \G \to G$. Denote by $n$ the order of
$\varphi(g)$ and write $k=\frac1n|G|$. Pick representatives $h_1,\dots,h_k$ for $G/\langle
g\rangle$. Then any $h\in G$ is of the form $h=h_ig^j$ for unique $i\in \{1,\dots,k\}$ and $j\in
\{0,\dots,n-1\}$.

Denoting the standard basis of $\Z^4$ again by $e_1,\dots,e_4$,
the linearity of $\lambda_G$ gives
 \[
 \lambda_G(h_jg^ke_i,h_{j'}g^{k'}e_{i'})=h_{j'}\lambda_G(g^ke_i,g^{k'}e_{i'})h_{j}^{-1}.
\]
Note that  by the definition of $\lambda_G$ we have $\lambda_G(g^ke_i,g^{k'}e_{i'})\in \Z[\langle
g\rangle]\in \Z[\G]$. Furthermore, $h_{j'}g^ih_j^{-1}=e$ if and only if $j'=j$ and $i=0$.
Therefore
$$
\lambda_G(h_jg^ke_i,h_{j'}g^{k'}e_{i'})_1=(h_{j'}\lambda_G(g^ke_i,g^{k'}e_{i'})h_{j}^{-1})_1=
\delta_{jj'}\,\lambda(e_ig^k,e_{i'}g^{k'})_1.
$$
This shows that the ordinary intersection form of $M_\varphi$ is the
direct sum of $k$ copies of $L_n$. This is clearly positive
definite. Taking $k$ copies of the vector $w'$
from the proof of Theorem \ref{thm:new} we get  a characteristic
vector of norm
$$
\hskip 1.5in (4n-8)k < 4|G|=\mbox{rank}(H_2(M_\varphi)). \hskip 1.5in \square
$$


\begin{thebibliography}{FQ}

\bibitem{CG} A. Casson and C. Gordon, {\em Cobordism of classical knots},
Progr. Math., 62,  A la recherche de la topologie perdue,  181--199, Birkh\"auser Boston, Boston,
MA, 1986.

\bibitem{C}  J.\,H. Conway,
{\em The Sensual (quadratic) Form},
Carrus Math. Monographs {\bf 26}, M.A.A. 1993.

\bibitem{CS}  J.\,H. Conway and N.\,J.\,A. Sloane,
{\em Sphere Packings, Lattices and Groups},\\
New York: Springer 1993.

\bibitem{D} S. Donaldson, {\em The orientation of Yang-Mills moduli spaces and $4$-manifold topology},\\
J. Differential Geom. {\bf 26} (1987), 397--428.

\bibitem{E} N. Elkies, {\em A characterization of the $\Z^n$ lattice}, Math. Res. Lett. {\bf 2} (1995), 321--326.

\bibitem{E'} \bysame, {\em Lattices and codes with long shadows}, Math. Res. Lett. {\bf 2} (1995), 643--651.

\bibitem{FT} R. Fintushel and R. Stern, {\em A fake {$4$}-manifold with {$\pi\sb 1={\bf Z}$} and {$b\sp+=4$}},\\ Turkish J. Math. {\bf 18} (1994), 1--6.

\bibitem{FQ}  M.\,H. Freedman and F. Quinn,
{\em Topology of 4-Manifolds},\\
Princeton Math. Series 39, Princeton, NJ, 1990.


\bibitem{F} S. Friedl, {\em Eta invariants as sliceness obstructions and their relation to Casson-Gordon
invariants},  Algebr. Geom. Topol.  4  (2004), 893--934.

\bibitem{G}  M. Gaulter, {\em Lattices without short characteristic vectors},
Math. Res. Lett. {\bf 5} (1998), 353--362.

\bibitem{HT}  I. Hambleton and P. Teichner, {\em A non-extended hermitian form over $\Z[\Z]$},\\
Manuscripta Mathematica {\bf 94} (1997),  435--442.

\bibitem{Ly} R. Lyndon, {\em Cohomology theory of groups with a single
defining relation},
\\ Ann. of Math. (2) {\bf 52}, (1950), 650--665.

\bibitem{MH}  J. Milnor and D. Husemoller,
{\em Symmetric Bilinear Forms},
New York: Springer 1973.

\bibitem{Q} H.-G. Quebbemann, {\em Definite lattices over real  algebraic function domains},\\
Math. Annalen {\bf 272} (1985), 461--475.

\end{thebibliography}
\end{document}